\documentclass[10pt]{article}
\usepackage{amssymb, amsmath, url, graphicx, setspace, geometry}
\usepackage{calrsfs, xcolor}
\usepackage{wasysym}

\def\3{\subset }
\def\4{\subseteq }
\def\<{\left<}
\def\>{\right>}

\def\bit{\begin{itemize}}
\def\eit{\end{itemize}}
\def\3{\subset }
\def\4{\subseteq }

\def\0{\leqno}

\def\barr{\begin{array}}
\def\earr{\end{array}}

\def\Z{{\rlap{$\kern2pt{\rm Z}$}{\rm Z}\,}}

\title{\bf Detecting solvability, supersolvability and CLT properties via an invariant based on conjugacy classes of subgroups}
\author{Mihai-Silviu Lazorec}

\begin{document}

\maketitle

\begin{abstract}
For a finite group $G$, denote by $k'(G)$ and $L(G)$ the number of conjugacy classes of subgroups and the subgroup lattice of $G$, respectively. Let $d'(G)=\frac{k'(G)}{|L(G)|}$ and $d^*(G)$ be the minimum value of $d'(S)$, when $S$ runs through all sections of $G$. In this paper we deduce some criteria on the nature of $G$. We show that if $d^*(G)>\frac{9}{59}$, then $G$ is solvable, while if $d^*(G)>\frac{1}{2}$, then $G$ is a supersolvable group. The last criterion is also valid when replacing ``supersolvable" with ``CLT". 
\end{abstract}



\section{Introduction}

All groups considered in this paper are finite. Given a group $G$, we denote by $k'(G)$ and $L(G)$ the number of conjugacy classes of subgroups and the subgroup lattice of $G$, respectively. The cyclic group of order $n$ and the dihedral group of order $2n$ are denoted by $C_n$ and $D_{2n}$, respectively.  
By a section of $G$, we mean a quotient group $H/N$, where $H, N\in L(G)$ and $N$ is a normal subgroup of $H$. We recall that $G$ is an Iwasawa group if all subgroups of $G$ are permutable or, equivalently, $G$ is a nilpotent modular group (see Exercise 3, p. 87 of \cite{12}). Also, $G$ is a CLT group if $G$ has a subgroup of order $d$ for any divisor $d$ of $|G|$. For a positive integer $n$, we denote by $\tau(n)$ and $\sigma(n)$ the number and the sum of all divisors of $n$, respectively.

The following ratio was studied in one of the author's previous papers (see \cite{8}):
$$d'(G)=\frac{k'(G)}{|L(G)|}.$$
It is obvious that $0<d'(G)\leq 1$. Also, the equality $d'(G)=1$ holds if and only if all subgroups of $G$ are normal. Hence $d'$ measures how close $G$ is from being a Dedekind group. It is known that the following sequence of inclusions between various classes of groups holds:
$$\text{Dedekind}\subset \text{Iwasawa}\subset \text{nilpotent}\subset \text{supersolvable}\subset \text{CLT}\subset \text{solvable}.$$

By the main result of \cite{2}, a group $G$ having only one conjugacy class of non-normal subgroups is isomorphic to one of the following:
\begin{itemize}
\item a modular $p$-group $M_{p^n}=\langle x, y \mid x^{p^{n-1}}=y^p=1,  yx=x^{p^{n-2}+1}y\rangle$ ($n\geq 4$ if $p=2$; $n\geq 3$ if $p$ is odd); 
\item a non-abelian split extension $G_{p, q, n}=H\rtimes K$ where $H\cong C_p$, $K\cong C_{q^{n-1}}, [H, \Phi(K)]=1$, $p, q$ are primes such that $q\mid p-1$ and $n\geq 2$.
\end{itemize} 
According to Proposition 2.2 and Corollary 2.3 of \cite{8}, one has
$$d'(G_{p,q,n})=\frac{2n}{2n+p-1}\xrightarrow[n \to \infty]{} 1.$$
This asymptotic behavior along with the fact that $G_{p, q, n}$ is non-nilpotent led to the impossibility of determining a constant $c\in (0, 1)$ such that if $d'(G)>c$, then $G$ is a nilpotent/Iwasawa/Dedekind group. This was a motivation for introducing the ratio 
$$d^*(G)=\min\{ d'(S) \mid S\text{ is a section of }G\},$$
which was used to obtain the following criteria on the nature of $G$ (see Theorems B, C, D in \cite{8}):
\begin{itemize}
\item If $d^*(G)>\frac{2}{3}=d^*(S_3)=d'(S_3)$, then $G$ is a nilpotent group;
\item If $d^*(G)>\frac{4}{5}=d^*(D_8)=d'(D_8)$, then $G$ is an Iwasawa group; 
\item Let $n\geq 3$ and $G$ be a $p$-group of order $p^n$. 
\begin{itemize}
\item Assume that $p=2$ and $n=3$. If $d^*(G)>\frac{4}{5}=d^*(D_8)=d'(D_8)$, then $G$ is a Dedekind group;
\item Assume that $n\geq 4$ if $p=2$, and $n\geq 3$ if $p$ is odd. If $d^*(G)>\frac{(n-2)(p+1)+4}{(n-1)(p+1)+2}=d^*(M_{p^n})=d'(M_{p^n})$, then $G$ is a Dedekind group. 
\end{itemize}
\end{itemize}
In each of the previous criteria, the outlined lower bound is the best possible one. 

In this paper, by making use of the classification of minimal simple groups, we are able to obtain a solvability criterion. Also, the classification of   minimal supersolvable groups with trivial Frattini subgroup leads us to establishing a condition which  guarantees that $G$ is a supersolvable group. More exactly, after recalling and justifying some preliminary results we prove the following statements in Section 2:\\

\textbf{Theorem [S].} \textit{Let $G$ be a group. If $d^*(G)>\frac{9}{59}=d^*(A_5)=d'(A_5)$, then $G$ is solvable.}\\

\textbf{Theorem [SS].} \textit{Let $G$ be a group. If $d^*(G)>\frac{1}{2}=d^*(A_4)=d'(A_4)$, then $G$ is supersolvable.}\\

Since any supersolvable group is a CLT group, an immediate consequence of Theorem [SS] is Corollary [CLT] below. Still, it is worth noting that a large part of the argument given in the proof of Theorem [SS] can be replicated to independently prove Corollary [CLT]. This separate proof is shared in Section 2 and it is based mainly on the fact that $G$ is a minimal non-CLT group if and only if it is a minimal non-supersolvable group. As we did not find it explicitly stated in literature, this equivalence is also justified in the same section.\\  

\textbf{Corollary [CLT].} \textit{Let $G$ be a group. If $d^*(G)>\frac{1}{2}=d^*(A_4)=d'(A_4)$, then $G$ is a CLT group.}\\

Since $A_5$ is non-solvable and $A_4$ is not supersolvable, nor CLT, the lower bounds $\frac{9}{59}$ and $\frac{1}{2}$ are the best possible ones. We note that the condition $d'(G)>\frac{9}{59}$ does not guarantee the solvability of $G$. By parsing through the groups $G$ with $|G|\leq 2000$ using GAP \cite{15}, we found 3 non-solvable groups, namely $SL(2, 5), C_7\times SL(2, 5), C_{11}\times SL(2, 5)$, each with a $d'$-value of $\frac{3}{19}$. Also, the inequality $d'(G)>\frac{1}{2}$ does not assure that $G$ is a supersolvable/CLT group. By searching through the groups $G$ with $|G|\leq 1000$, one can find 34 non-supersolvable groups (33 of them being non-CLT) with $d'$-values falling within the interval $[\frac{40}{79}, \frac{19}{33}]$. The following table lists these groups along with their $d'$-values in increasing order.

\begin{center}
\noindent\begin{tabular}{ |p{7cm}|p{2cm}|p{2cm}|}
  \hline
IdGroup($G$) in GAP's Small Groups library & CLT group & $d'(G)$\\
 \hline
[432, 41] & Yes & $40/79$\\
\hline
[972, 832] & No & $251/495$\\
\cline{1-1}\cline{3-3}
[972, 147] &  & $23/45$\\
\cline{1-1}\cline{3-3}
[324, 133]  &  & $64/125$\\
\cline{1-1}\cline{3-3}
[108, 20], [540, 59], [ 756, 109 ] &  & $13/25$\\
\cline{1-1}\cline{3-3}
[72, 3], [360, 14], [504, 22], [792, 12], [936, 25] &  & $11/21$\\
\cline{1-1}\cline{3-3}
[972, 118] &  & $19/36$\\
\cline{1-1}\cline{3-3}
[36, 3], [180, 8], [252, 10], [396, 6], [468, 13], [612, 8], [684, 10], [828, 6], [900, 8], [900, 65],  &  & $8/15$\\
\cline{1-1}\cline{3-3}
[972, 539] &  & $97/180$\\
\cline{1-1}\cline{3-3}
[648, 88] &  & $13/24$\\
\cline{1-1}\cline{3-3}
[324, 44] &  & $19/35$\\
\cline{1-1}\cline{3-3}
[108, 3], [540, 8], [756, 10] &  & $11/20$\\
\cline{1-1}\cline{3-3}
[216, 3], [972, 146] &  & $5/9$\\
\cline{1-1}\cline{3-3}
[324, 3] &  & $14/25$\\
\cline{1-1}\cline{3-3}
[972, 3] &  & $17/30$\\
\cline{1-1}\cline{3-3}
[648, 3] &  & $19/33$\\
\hline
 \end{tabular}
 \end{center}
Based on these facts, we end this section by posing the following question.\\

\textbf{Open problem.} \textit{Can one find a constant $c\in (0, 1)$ such that if $d'(G)>c$, then $G$ is a solvable/supersolvable/CLT group?}

\section{Proofs of the main results}

To prove Theorem [S], we first list 4 preliminary results. The first one (see Corollary 1 in \cite{16}) outlines the classification of minimal simple groups (i.e. non-abelian simple groups all of whose proper subgroups are solvable). The following two rely on the classifications of the maximal subgroups of projective special linear groups (see Exercise 7, p. 417 of \cite{14}) and of Suzuki groups (see Theorem 4.1 in \cite{17}). Based on Theorems 3.1 and 3.3 of \cite{3}, the last lemma establishes some results on the $d'$-value of a dihedral group. The proof of Theorem [S] is given afterwards.\\

\textbf{Lemma 2.1.} \textit{Let $G$ be a minimal simple group. Then $G$ is isomorphic to one of the following groups:
\begin{itemize}
\item $PSL(2, 2^t)$, where $t$ is a prime;
\item $PSL(2, 3^t)$, where $t\geq 3$ is a prime;
\item $PSL(2, p)$, where $p\geq 7$ is a prime such that $p^2+1\equiv 0 \ (mod \ 5)$;
\item $Sz(2^t)$, where $t\geq 3$ is a prime;
\item $PSL(3,3)$.
\end{itemize}}

\textbf{Lemma 2.2.} \textit{Let $q=p^t$, where $p$ is a prime and $t\geq 1$.
\begin{itemize}
\item If $q\geq 4$ is even, then $PSL(2, q)$ has a maximal subgroup isomorphic to $D_{2(q-1)}$;
\item If $q\geq 13$ is odd, then $PSL(2, q)$ has a maximal subgroup isomorphic to $D_{q-1}$.
\end{itemize}}

\textbf{Lemma 2.3.} \textit{Let $q=2^t$, where $t\geq 3$ is odd. Then $Sz(q)$ has a maximal subgroup isomorphic to $D_{2(q-1)}$.}\\

\textbf{Lemma 2.4.} \textit{Let $n\geq 3$. Then
$$d'(D_{2n})=\frac{4\tau(n)-\tau(2n)}{\tau(n)+\sigma(n)} \text{ \ and \ } d'(D_{2n})< \frac{12\sqrt[3]{n}}{n+3}.$$}

\textbf{Proof.} It is known that 
$$D_{2n}=\langle x, y\mid x^n=y^2=1, yx=x^{-1}y\rangle.$$
According to Theorem 3.1 of \cite{3}, for each divisor $d$ of $n$, $D_{2n}$ has: 
\begin{itemize}
\item[--] one subgroup of index $2d$: $\langle x^d \rangle$;
\item[--] $d$ subgroups of index $d$: $\langle x^d, x^iy\rangle$, with $i\in\{ 0, 1,\ldots, d-1\}.$
\end{itemize}
Hence, we have $$|L(D_{2n})|=\tau(n)+\sigma(n).$$ 

By Theorem 3.3 of \cite{3}, if $n$ is odd, then all subgroups of a fixed index form a single conjugacy class, yielding $$k'(D_{2n})=\tau(2n).$$ 
By the same result, if $n$ is even and $e\mid 2n$, there is a unique conjugacy class of subgroups of index $e$ unless $e$ is even and $e\mid n$, in which case there are 3 such classes. Hence,
\begin{align*}
k'(D_{2n})=\sum\limits_{\substack{e\mid 2n\\ e\equiv 1 \ (mod \ 2)}}1+\sum\limits_{\substack{e\mid 2n\\ e\nmid n\\e\equiv 0 \ (mod \ 2)}}1+\sum\limits_{\substack{e\mid 2n\\ e\mid n\\e\equiv 0 \ (mod \ 2)}}3=\tau(2n)+\sum\limits_{\substack{ e\mid n\\e\equiv 0 \ (mod \ 2)}}2=\tau(2n)+2\tau(\frac{n}{2}).
\end{align*}

By the standard multiplicativity property of the $\tau$ function, it is easy to check that 
$$k'(D_{2n})=4\tau(n)-\tau(2n) \text{ \ and \ } k'(D_{2n})<3\tau(n)$$
regardless of the parity of $n$.  
It is clear that $\tau(n)+\sigma(n)\geq n+3$. Also, it is known that $\tau(n)<4\sqrt[3]{n}$ (see, for instance, Exercise 3.5.1 in \cite{13}). From the above facts, we deduce that
$$d'(D_{2n})=\frac{4\tau(n)-\tau(2n)}{\tau(n)+\sigma(n)}<\frac{3\tau(n)}{n+3}<\frac{12\sqrt[3]{n}}{n+3},$$
as desired.
\hfill\rule{1,5mm}{1,5mm}\\


\textbf{Proof of Theorem [S].} For the sake of contradiction, let $G$ be a non-solvable group of minimal order such that $d^*(G)>\frac{9}{59}$. Let $S\not\cong G$ be a section of $G$. It is easy to check that any section of $S$ is isomorphic to a section of $G$. Hence, $d^*(S)\geq d^*(G)$, which implies that  
\begin{equation}\label{rel1}
d^*(S)>\frac{9}{59}.
\end{equation}
It follows that $S$ is solvable since $|S|<|G|$. In particular, any proper subgroup $H$ of $G$ is solvable.

Let $N$ be a maximal normal subgroup of $G$. Then $G/N$ is simple and it cannot be abelian since this would mean that both $N$ and $G/N$ are solvable and, consequently, $G$ would be solvable, a contradiction. The proper subgroups of $G/N$ are $H/N$, where $H$ is a proper subgroup of $G$ containing $N$. Since $H$ is solvable, so is $H/N$. These facts lead us to stating that $G/N$ is a minimal simple group. Thus, $G/N$ must be isomorphic to one of the groups listed in Lemma 2.1. To finish the proof, it suffices to show that $d^*(G/N)\leq\frac{9}{59}$ in all cases, contradicting \eqref{rel1}.

$\bullet$ $G/N\cong PSL(2, 4)$ or $G/N\cong PSL(3,3)$

One may use GAP to obtain
$$d^*(PSL(2,4))=\frac{9}{59} \hspace{2cm} d^*(PSL(3,3))=\frac{51}{6374}<\frac{9}{59}.$$

$\bullet$ $G/N\cong PSL(2, q)$ or $G/N\cong Sz(q)$, where $q=2^t$ and $t\geq 3$ is a prime

For $q\in \{2^3, 2^5, 2^7\}$, we summarize the values of $d^*(G/N)$ is the following table:
\begin{center}
\noindent\begin{tabular}{ |p{2cm}|p{4cm}|p{4cm}|}
  \hline
$q$ & $d^*(PSL(2, q))$ & $d^*(Sz(q))$\\
\hline
$8$ & $6/193$ &  $22/17295$\\
\hline
$32$ & $3/2791$ & $132/21170191$\\
\hline
$128$ & $121/2095631$ & $\leq 2/65=d'(D_{254})<9/59$\\
\hline
 \end{tabular}
 \end{center}
We remark that all values are less than $\frac{9}{59}$, as desired. 

In what follows, assume that $q\geq 2^{11}$. According to Lemmas 2.2 and 2.3, $G/N$ has a subgroup isomorphic to $D_{2(q-1)}$. Hence, $d^*(G/N)\leq d'(D_{2(q-1)})$. By Lemma 2.4, we get
\begin{align}\label{rel2}
d'(D_{2(q-1)})<\frac{12\sqrt[3]{q-1}}{q+2}.
\end{align}
Let $f:[2^{11}, \infty)\longrightarrow \mathbb{R}$ be a function given by $$f(x)=\frac{\sqrt[3]{x-1}}{x+2}, \ \forall \ x\in [2^{11}, \infty).$$
We deduce that $f$ is strictly decreasing since
$$f'(x)=\frac{5-2x}{3\sqrt[3]{(x-1)^2}(x+2)^2}<0, \ \forall \ x\in (2^{11}, \infty).$$
It follows from \eqref{rel2} that
$$d'(D_{2(q-1)})<\frac{12\sqrt[3]{2^{11}-1}}{2^{11}+2}<\frac{9}{59},$$ and, consequently,
$$d^*(G/N)<\frac{9}{59}.$$

$\bullet$ $G/N\cong PSL(2, q)$, where $q=3^t$ for a prime $t\geq 3$, or $q=p$ for a prime $p\geq 7$ satisfying $p^2+1\equiv 0 \ (mod \ 5)$ 
  
For $q=7$, we have
$$d^*(PSL(2,7))=\frac{15}{179}<\frac{9}{59}.$$

Assume that $q\geq 13$, so $PSL(2, q)$ has a subgroup isomorphic to $D_{q-1}$ according to Lemma 2.2. 
If $q\geq 487$, by an argument similar to the one conducted in the previous case, we have 
$$d'(D_{q-1})<\frac{24\sqrt[3]{q-1}}{\sqrt[3]{2}(q+5)}\leq \frac{24\sqrt[3]{486}}{492\sqrt[3]{2}}<\frac{9}{59},$$
implying
$$d^*(PSL(2, q))<\frac{9}{59}.$$

If $q=3^t\in [13, 487)$ for a prime $t\geq 3$, we have
$$d^*(PSL(2, 27))=\frac{8}{2643}<\frac{9}{59} \hspace{2cm} d^*(PSL(2, 243))=\frac{13}{812482}<\frac{9}{59}.$$

For the 45 possible choices of $q=p\in [13, 487)$ being a prime that satisfies $p^2+1 \equiv 0 \ (mod \ 5)$, we share our results in the following table:
\footnote{$d^*(PSL(2, q))$ was computed only if $d'(D_{q-1})>\frac{9}{59}$; otherwise, a ``$\star$" was placed in the corresponding cell of the table meaning that the inequality $d^*(PSL(2, q))< \frac{9}{59}$ holds as a consequence of $d'(D_{q-1})<\frac{9}{59}$.} 
 \begin{center}
\noindent\begin{tabular}{ |p{0.4cm}|p{1.2cm}|p{2.1cm}|p{0.4cm}|p{1.2cm}|p{2.1cm}|p{0.4cm}|p{1.2cm}|p{2.1cm}|}
\hline
$q$ & $d'(D_{q-1})$ & $d^*(PSL(2, q))$ & $q$ & $d'(D_{q-1})$ & $d^*(PSL(2, q))$ & $q$ & $d'(D_{q-1})$ & $d^*(PSL(2, q))$\\
\hline
13 & 5/8 & 8/471 & 137 & 4/33 & $\star$ & 307 & 1/20 & $\star$ \\
\hline
17 & 11/19 & 1/110 & 157 & 5/44 & $\star$ & 313 & 8/101 & $\star$ \\
\hline
23 & 2/7 & 23/5915 & 163 & 5/63 & $\star$ & 317 & 5/122 & $\star$ \\
\hline
37 & 1/3 & 23/17731 & 167 & 2/43 & $\star$ & 337 & 11/124 & $\star$ \\
\hline
43 & 2/9 & 10/12731 & 173 & 5/68 & $\star$ & 347 & 1/44 & $\star$ \\
\hline
47 & 2/13 & 29/48837 & 193 & 17/132 & $\star$ & 353 & 14/191 & $\star$ \\
\hline
53 & 5/23 & 10/21627 & 197 & 5/59 & $\star$ & 367 & 2/63 & $\star$ \\
\hline
67 & 2/13 & 10/39801 & 223 & 2/39 & $\star$ & 373 & 5/98 & $\star$ \\
\hline
73 & 6/25 & 38/176087 & 227 & 1/29 & $\star$ & 383 & 2/97 & $\star$ \\
\hline
83 & 1/11 & $\star$ & 233 & 2/27 & $\star$ & 397 & 1/16 & $\star$ \\
\hline
97 & 14/67 & 45/451547 & 257 & 23/263 & $\star$ & 433 & 1/14 & $\star$ \\
\hline
103 & 2/19 & $\star$ & 263 & 2/67 & $\star$ & 443 & 1/32 & $\star$ \\
\hline
107 & 1/14 & $\star$ & 277 & 5/74 & $\star$ & 457 & 8/143 & $\star$ \\
\hline
113 & 11/64 & 39/622753 & 283 & 2/49 & $\star$ & 463 & 2/49 & $\star$ \\
\hline
127 & 6/55 & $\star$ & 293 & 5/113 & $\star$ & 467 & 1/59 & $\star$ \\
\hline
\end{tabular}
\end{center}
Based on our findings, we conclude that in all these 45 cases the inequality
$$d^*(PSL(2, q))<\frac{9}{59}$$
holds, as desired.
\hfill\rule{1,5mm}{1,5mm}\\

Our following objective is to prove Theorem [SS]. First, we recall a classification of the minimal non-supersolvable groups with trivial Frattini subgroup (see Theorem 1 in \cite{10} and Lemma 2.2 in \cite{6}). For a full classification of the minimal non-supersolvable groups, one may check Theorems 9, 10 in \cite{1}. Finally, the sixth lemma gathers auxiliary  properties that will be useful in proving Theorem [SS].\\ 

\textbf{Lemma 2.5} \textit{Let $p, q, r$ be distinct primes and $G$ be a minimal non-supersolvable group with $\Phi(G)=\{1\}$. Then $G$ is isomorphic to one of the following:
\begin{itemize}
\item[(1)] a Frobenius group $C_p^k\rtimes C_q$, where $k\geq 2$ is the smallest integer such that $p^k\equiv 1 \ (mod \ q)$;
\item[(2)] a Frobenius group $C_p^q\rtimes C_{q^{n+1}}$, where $n\geq 1, p\equiv 1 \ (mod \ q^n), p\not\equiv 1 \ (mod \ q^{n+1})$;
\item[(3)] a Frobenius group $C_p^2\rtimes Q_8$, where $p\equiv 1 \ (mod \ 4)$;
\item[(4)] $C_p^q\rtimes (C_{q^n}\rtimes C_q)$, where $n\geq 2, p\equiv 1 \ (mod \ q^n)$;
\item[(5)] $C_p^q\rtimes ((C_q\times C_q)\rtimes C_q)$, where $p\equiv 1 \ (mod \ q)$ and $q\geq 3$;
\item[(6)] $C_p^r\rtimes (C_q\rtimes C_{r^n})$, where $n\geq 1, p\equiv 1 \ (mod \ qr^n), q\equiv 1 \ (mod \ r)$.
\end{itemize}
In all cases, the Sylow $p$-subgroup of $G$ is the Fitting subgroup $F(G)$.}\\

\textbf{Lemma 2.6.} \textit{Let $G$ be a group.
\begin{itemize}
\item[i)] If $G/\Phi(G)$ is supersolvable, then $G$ is also supersolvable;
\item[ii)] If $G$ is solvable, then $C_G(F(G))\leq F(G)$;
\item[iii)] If $H$ and $K$ are permutable subgroups of $G$, then $[H, K]$ is a normal subgroup of $HK$;
\item[iv)] If $G$ is supersolvable, then it has a cyclic normal subgroup of prime order;  
\item[v)] Let $A$ be a group acting via automorphisms on $G$. If $(|G|, |A|)=1$ and $G$ is abelian, then $G=[G, A]\times C_G(A)$;
\item[vi)] Let $p, q$ be primes such that $p>q$ and $|G|=p^t q$, where $t\geq 1$. If $G$ is a non-nilpotent supersolvable group and its Sylow $p$-subgroup $P$ is abelian, then $G$ has a subgroup isomorphic to the non-abelian group of order $pq$.
\end{itemize}}

\textbf{Proof.} Items \textit{i)} and \textit{ii)} are parts of Theorems 8.6 and 4.2 of \cite{7} (see p. 754, p. 292). The third statement follows by Lemmas 3.11 and 1.6 of the same reference (see p. 18, 268). The fourth and fifth results are consequences of 5.4.8 in \cite{11} and Theorem 2.3, p. 177 in \cite{5}, respectively.

Regarding item \textit{vi)}, it is clear that the result holds for $t=1$, so we may assume that $t\geq 2$. Let   $Q\cong C_q$ be a Sylow $q$-subgroup of $G$. We have $G=PQ$ and, since $G$ is non-nilpotent, it follows that $[P, Q]\neq \{ 1\}$. Moreover, $[P, Q]$ is a normal subgroup of $G$ according to \textit{iii)}. Since $[P, Q]$ is a non-trivial normal subgroup of $G$, it must contain  a minimal normal subgroup $N$ of $G$. By \textit{iv)}, $N$ must be cyclic of prime order. Since $G$ is non-nilpotent, we have $N\cong C_p$. Then $NQ$ is a subgroup of order $pq$ of $G$. It remains to show that $NQ$ is non-abelian. 

By \textit{v)}, we have $P=[P, Q]\times C_P(Q)$, so $[P, Q]\cap C_P(Q)=\{ 1\}$. Since $N$ is a subgroup of $[P, Q]$, it follows that $N\cap C_P(Q)=\{ 1\}$. As $N$ is contained in $P$, we conclude that $NQ$ is non-abelian.
\hfill\rule{1,5mm}{1,5mm}\\

\textbf{Proof of Theorem [SS].} Assuming the contrary, let $G$ be a non-supersolvable group of minimal order with $d^*(G)>\frac{1}{2}$. As in the proof of Theorem [S], one can show that $d^*(S)>\frac{1}{2}$
for any section $S\not\cong G$ of $G$. By the minimality of $|G|$, any such section and, in particular, any proper subgroup of $G$ would be supersolvable, making $G$ a minimal non-supersolvable group. Based on these facts, if $\Phi(G)\neq \{1\}$, then $G/\Phi(G)$ is supersolvable. By Lemma 2.6 \textit{i)}, $G$ would also be supersolvable, a contradiction. Hence, $\Phi(G)=\{1\}$, so $G$ is isomorphic to one of the groups listed in Lemma 2.5. To complete the proof, it suffices to show that $d^*(G)\leq\frac{1}{2}$ in each case. We denote the Sylow $p$-subgroup of $G$ by $P$.

$\bullet$ $G$ is of type \textit{(3), (4), (5)} or \textit{(6)}

Based on the arithmetical conditions associated with these 4 isomorphism classes, we remark that $p\geq 5$. Let $Q\cong C_q$ be a subgroup of $G$ ($q=2$ if $G$ is of type \textit{(3)}). Based on the structure of $G$, it is easy to see that $H=PQ$ is a proper subgroup of $G$, so it is supersolvable. 

Assume that $H$ is nilpotent. Then $H$ is abelian, so $Q\subset C_G(P)$. Since $d^*(G)>\frac{1}{2}>\frac{9}{59}$, $G$ is solvable by Theorem [S] (this is also guaranteed by the fact that any minimal non-supersolvable group is solvable; see Theorem 12 in \cite{1} or Hilfssatz C in \cite{4}). By Lemma 2.5, we know that $F(G)=P$. According to Lemma 2.6 \textit{ii)}, it follows that $C_G(P)\subset P$, so $Q\subset P$, a contradiction. 

Hence, $H$ is a non-nilpotent supersolvable group and its Sylow $p$-subgroup $P$ is abelian. By Lemma 2.6 \textit{vi)}, $H$ has a subgroup $K$ isomorphic to the non-abelian group of order $pq$. Then,
$$d^*(G)\leq d'(K)=\frac{4}{p+3}\leq\frac{4}{8}=\frac{1}{2}.$$

$\bullet$ $G$ is of type \textit{(2)}

If $n\geq 2$, the argument given in the previous case can be repeated. The same is valid if $n=1$ and $p\geq 5$. Hence, the only left possibility is that $(p, q, n)=(3, 2, 1)$. Then $G\cong C_3^2\rtimes C_4$ (SmallGroup(36, 9)) and $$d^*(G)=\frac{5}{19}<\frac{1}{2}.$$ 

$\bullet$ $G$ is of type \textit{(1)}

Let $Q\cong C_q$ be a Frobenius complement of $G$. Then 
\begin{align}\label{rel3}
C_Q(x)=\{1\}, \ \forall \ x\in P\setminus\{1\}.
\end{align}
First, we show that $G$ does not have a subgroup of order $p^iq$ for any $i\in \{1, 2, \ldots, k-1\}$. Assume the contrary and let $H_i$ be a subgroup such that $|H_i|=p^iq$. Since $k$ is the smallest integer such that $p^k\equiv 1 \ (mod \ q)$, without loss of generality, we may assume that $Q$ is the Sylow $q$-subgroup of $H_i$. Let $P_i=H_i\cap P$. Since $P$ is an abelian normal subgroup of $G$, it is easy to check that $P_i$ is an abelian normal subgroup of  $H_i$. Then,
$$H_i/P_i=H_i/(H_i\cap P)\cong H_iP/P\leq G/P\cong C_q.$$
Since $H_i\not\subset P$, we get $|H_i/P_i|=q$ and, consequently, $|P_i|=p^i$. It follows that $P_i$ is the Sylow $p$-subgroup of $H_i$, so $H_i\cong P_i\times Q$. Since $P_i\subset P$, there is $x\in P\setminus \{1\}$ such that $C_Q(x)=Q$, contradicting \eqref{rel3}.

Therefore, 
$$L(G)=L(P)\cup \{Q^g \mid g\in G\}\cup \{G\}.$$
The size of the conjugacy class of $Q$ is $p^k$.
By the proof of Theorem B in \cite{8}, we have
$$d'(G)=\frac{a_{p, k}+4q-2}{q(a_{p, k}+p^k+1)},$$
where $a_{p, k}=|L(P)|$. 
Then $$d'(G)\leq\frac{1}{2}\Longleftrightarrow (q-2)a_{p, k}+(p^k-7)q+4\geq 0.$$
Since $q\geq 2$, the last inequality clearly holds if $p\geq 3$. The same is valid if $p=2$ and $k\geq 3$. The only left possibility is $(p, q, k)=(2, 3, 2)$. Then $G\cong A_4$ and $d'(G)=\frac{1}{2}$. 

Therefore, if $G$ is of type (\textit{1}), we have
$$d^*(G)\leq\frac{1}{2},$$
and the proof is complete. 
\hfill\rule{1,5mm}{1,5mm}\\

We end our paper by showing that Corollary [CLT] may be proved even without having Theorem [SS] at our disposal. To achieve this, we first justify that the classes of minimal non-supersolvable groups and minimal non-CLT groups coincide. This result is obtained by making use of the following lemma (see Theorem 1 in \cite{9}).\\

\textbf{Lemma 2.7.} \textit{Let $G$ be a group. Then $G$ is supersolvable if and only if all subgroups of $G$ are CLT groups.}\\

\textbf{Lemma 2.8} \textit{Let $G$ be a group. Then $G$ is minimal non-supersolvable if and only if $G$ is a minimal non-CLT group.}\\

\textbf{Proof.} Suppose that $G$ is minimal non-supersolvable. Then all its proper subgroups are supersolvable, so they are also CLT-groups. By Lemma 2.7, since $G$ is non-supersolvable, it has a non-CLT subgroup $H$. It follows that the only possibility is $H=G$, so $G$ is a non-CLT group all of whose proper subgroups are CLT.

Conversely, let $G$ be a minimal non-CLT group. Then $G$ is a non-CLT group and, consequently, it is non-supersolvable. Let $H$ be a proper subgroup of $G$. By the hypothesis, all subgroups of $H$ are CLT. Then $H$ is supersolvable according to Lemma 2.7. It follows that $G$ is a minimal non-supersolvable group.
\hfill\rule{1,5mm}{1,5mm}\\

\textbf{An alternative proof of Corollary [CLT].} Assuming the contrary, let $G$ be a non-CLT group of minimal order such that $d^*(G)>\frac{1}{2}$. Analogous to the proof of Theorem [S], it follows that any section $S\not\cong G$ of $G$ is CLT. It follows that $G$ is a minimal non-CLT group. Hence, $G$ is minimal non-supersolvable by Lemma 2.8. Also, if we assume that $\Phi(G)\neq \{1\}$, then all subgroups of $G/\Phi(G)$ are CLT. Thus, $G/\Phi(G)$ is supersolvable by Lemma 2.7. According to Lemma 2.6 \textit{i)}, it follows that $G$ is supersolvable, a contradiction.

We conclude that $G$ is a minimal non-supersolvable group with trivial Frattini subgroup, so it is isomorphic to one of the groups listed in Lemma 2.5. From this point on, the argument proceeds exactly as in the proof of Theorem [SS].   
\hfill\rule{1,5mm}{1,5mm}\\

\bigskip\noindent {\bf Declarations}

\bigskip\noindent {\bf  Funding.} The author did not receive support from any organization for the submitted work.

\bigskip\noindent {\bf Conflicts of  interests.} The author declares that there is  no conflict of interest.

\bigskip\noindent {\bf  Data availability statement.} The manuscript has no associated data.

\vspace*{3ex}
\small
\hfill
\begin{minipage}[t]{7cm}
Mihai-Silviu Lazorec \\
Faculty of  Mathematics \\
"Al. I. Cuza" University \\
Ia\c si, Romania \\
e-mail: {\tt silviu.lazorec@uaic.ro}
\end{minipage}
\end{document}